\documentclass{amsart}
\newtheorem{thm}{Theorem}
\newtheorem{cor}[thm]{Corollary}
\newtheorem{lem}[thm]{Lemma}

\newtheorem{prop}[thm]{Proposition}

\newtheorem{rem}[thm]{Remark}
\newtheorem{ex}[thm]{Example}

\newtheorem{defn}[thm]{Definition}

\numberwithin{equation}{section}

\begin{document}

\title[Inf-convolution and convex functions on Riemannian manifolds]{
Inf-convolution and regularization of convex functions on
Riemannian manifolds of nonpositive curvature}
\author{Daniel Azagra and Juan Ferrera}

\address{Departamento de An\'{a}lisis Matem\'{a}tico\\ Facultad de
Matem\'{a}ticas\\ Universidad Complutense\\ 28040 Madrid, Spain}

\date{May 17, 2005}

\email{azagra@mat.ucm.es, ferrera@mat.ucm.es}


\subjclass[2000]{}

\begin{abstract}
We show how an operation of  inf-convolution can be used to
approximate convex functions with $C^{1}$ smooth convex functions
on Riemannian manifolds with nonpositive curvature (in a manner
that not only is explicit but also preserves some other properties
of the original functions, such as ordering, symmetries, infima
and sets of minimizers), and we give some applications.
\end{abstract}

\maketitle

\section{Introduction and main results}

Smooth approximation is an old subject. Its importance lies on the
fact that most analytical tools for studying the properties of
functions defined on a normed space or on a Riemannian manifold
require some degree of differentiability of the considered
functions. However, many functions which arise naturally from
geometrical problems on manifolds are only continuous (or even
merely lower semicontinuous). One is thus tempted to approximate
those functions by smooth functions to which one can apply more
powerful analytical methods and obtain some information about the
behavior of the approximations which will hopefully be shared with
the original (nonsmooth) functions.

The theory of convex functions is also an old subject which plays
an important role in many areas of mathematics. In Riemannian
geometry it has been used, for instance, in the investigation of
the structure of noncompact manifolds of positive curvature by
Cheeger, Gromoll, Greene, Meyer, Siohama, Wu and others, see
\cite{GromollMeyer, CheegerGromoll, Greene1, Greene2, Greene3,
Greene4} (recall that a function $f:M\to\mathbb{R}$ defined on a
Riemannian manifold $M$ is said to be {\em convex} provided the
function $\mathbb{R}\ni t\mapsto f\circ\gamma(t)\in \mathbb{R}$ is
convex for every geodesic $\gamma$ on $M$). The existence of
global convex functions on a Riemannian manifold has strong
geometrical and topological implications; for instance it is shown
in \cite{Greene1} that every two-dimensional manifold which admits
a global convex function which is locally nonconstant must be
diffeomorphic to the plane, the cylinder, or the open M\"{o}bius
strip.

Along with the papers cited above and the references therein, we
must mention the important work of Bangert's on convex sets and
convex functions on Riemannian manifolds, see \cite{Bangert,
Bangert2, Bangert3}; he showed in particular that Alexandroff's
Theorem about almost everywhere twice differentiability of convex
functions on $\mathbb{R}^{n}$ can be extended to convex functions
on finite-dimensional Riemannian manifolds (providing as well a
smart proof of Alexandroff's theorem on $\mathbb{R}^{n}$).

The aim of the present paper is to study to what extent one of the
most useful methods for regularizing convex functions on normed
spaces, namely that of infimal convolution, can be successfully
used in the setting of Riemannian manifolds.

Let us first have a quick look at the three main methods (that is,
partitions of unity, integral convolution with a sequence of
mollifiers, and inf and sup convolution formulae) that are used in
normed spaces to approximate continuous functions by smooth
functions, and see how they can be adapted to the case when one
wants to regularize a convex function $f$ defined on a Riemannian
manifold $M$.

Partitions of unity are useless in this setting because, even in
the case when $M=\mathbb{R}^{n}$, they do not preserve convexity
of the function $f$.

The integral convolution with a sequence of mollifiers reveals
itself as the perfect tool when $M=\mathbb{R}^{n}$, because the
integral convolution of a convex function $f$ with any integrable
function $g$ with compact support, that is,
    $$
    f*g(x)=\int_{\mathbb{R}^{n}}f(x-y)g(y)dy,
    $$
is a convex function. In \cite{Greene5, Greene3, Greene4} Greene
and Wu studied to what extent those integral convolutions with
mollifiers can be used to regularize convex functions defined on
finite-dimensional Riemannian manifolds $M$ (and applied this
method to prove several theorems about the structure of complete
noncompact manifolds of positive curvature). It turns out that
this method works out in Riemannian manifolds only when the
original function $f$ is {\em strictly convex}. More precisely,
let $\kappa:\mathbb{R}\to\mathbb{R}$ be a nonnegative $C^{\infty}$
function with support in $[-1,1]$, constant on a neighborhood of
$0$ and satisfying $\int_{\mathbb{R}^{n}}\kappa(\|x\|)dx=1$, and
 let us define the functions
    $$
    \varphi_{\varepsilon}(p)=
    \frac{1}{\varepsilon^{n}}\int_{v\in TM_{p}}f(\exp_{p}(v))\,
    \kappa\left(\frac{\|v\|_{p}}{\varepsilon}\right)d\mu_{p},
    $$
where $d\mu_{p}$ is the measure on the tangent space $TM_{p}$
obtained from the Riemannian metric of $M$. Greene and Wu showed
that if $f:M\to\mathbb{R}$ is a convex function defined on an
$n$-dimensional Riemannian manifold $M$ and $K$ is a compact
subset of $M$, then there exists an open neighborhood of $K$ and
an $\varepsilon_{0}>0$ such that the functions
$\varphi_{\varepsilon}:U\to\mathbb{R}$ defined above are
$C^\infty$ smooth, converge to $f$ uniformly on $K$ as
$\varepsilon\to 0$, and are {\em approximately convex} in the
sense that
    $$
    \liminf_{\varepsilon\to 0}\big(\inf
    \frac{d^{2}}{dt^{2}}\varphi_{\varepsilon}(\gamma(t))\Big|_{t=0}
    \big)\geq 0,
    $$
where the infimum is taken over all geodesics $\gamma(t)$
parameterized by arc length and with $\gamma(0)\in K$.

Now, a $C^{\infty}$ function $\varphi:M\to\mathbb{R}$ is called
{\em strictly convex} if its second derivative along any geodesic
is strictly positive everywhere on the geodesic. A (not
necessarily smooth) function $f:M\to\mathbb{R}$ is then said to be
{\em strictly convex} provided that for every $p\in M$ and every
$C^{\infty}$ strictly convex function $\varphi$ defined on a
neighborhood of $p$ there is some $\varepsilon>0$ such that
$f-\varepsilon\varphi$ is convex on the neighborhood. With this
terminology, the above result implies that if $f:M\to\mathbb{R}$
is strictly convex then for every compact subset $K$ of $M$ there
exists a sequence of strictly convex $C^{\infty}$ functions
$\varphi_{n}=\varphi_{\varepsilon_{n}}$ such that
$f=\lim_{n\to\infty}\varphi_{n}$ uniformly on $K$.

However, as Greene and Wu pointed out, this method cannot be used
when $f$ is not strictly convex, and the problem of approximating
(not necessarily strictly) convex functions by smooth convex
functions on Riemannian manifolds is open. That is one main
limitation of the integral convolution technique on Riemannian
manifolds. Another drawback of this method is the fact that it
does not apply to functions defined on infinite-dimensional
manifolds (even in the case when $M$ is the Hilbert space).

We are left with the third method: infimal convolution. It is well
known that if $f:X\to\mathbb{R}\cup\{+\infty\}$ is a lower
semicontinuous convex function defined in $X=\mathbb{R}^{n}$ or in
any infinite-dimensional reflexive Banach space $X$ (such as, for
instance, the separable Hilbert space), then the inf-convolution
formula
    $$
    f_{k}(x):=\inf_{y\in X}\{f(y)+k\|x-u\|^{2}\},
    $$
where $\|\cdot\|$ is any equivalent differentiable norm in $X$,
defines a sequence of $C^1$ smooth convex functions $f_{k}$ which
converge to $f$ as $k\to+\infty$ (uniformly on bounded sets if $f$
is uniformly continuous on bounded sets).  In fact, a clever
combination of inf and sup convolutions allows to show that if $f$
is a (not necessarily convex) function which is uniformly
continuous on bounded sets of a superreflexive Banach space then
$f$ can be approximated by $C^1$ smooth functions with uniformly
continuous derivatives uniformly on bounded sets; this was shown
first by Lasry and Lions \cite{Lasry-Lions} in the case when $X$
is the Hilbert space and then by Cepedello-Boiso \cite{C1, C2} for
any superreflexive Banach space $X$.

In this situation it is natural to ask whether infimal convolution
formulae can be used to regularize convex functions defined on
Riemannian manifolds (either finite or infinite-dimensional). That
is the question we try to address in this paper. Let us describe
the main results that we will show in the following sections. If
$M$ is a complete Riemannian manifold and $d$ is the geodesic
distance on $M$, for any function
$f:M\to\mathbb{R}\cup\{+\infty\}$, and for $\lambda>0$ we define
the function
$$
f_{\lambda}(x)=\inf_{y\in M}\{f(y)+\frac{1}{2\lambda}d(x,y)^{2}\}
$$
for every $x\in M$. In Section 2 we collect some general
properties of $f_{\lambda}$ that do not depend on the geometry of
$M$ and that we will need to use later in our proofs. We show for
instance how the inf defining $f_{\lambda}(x)$ can be restricted
to a suitable ball $B(x,r_{x})$, and then use the estimates on the
radius $r_{x}$ to show that $\lim_{\lambda\to
0^{+}}f_{\lambda}(x)=f(x)$ pointwise whenever $f(x)<+\infty$, and
that if $f$ is uniformly continuous on bounded sets then
$f_{\lambda}$ converges to $f$ uniformly on bounded sets. We also
see that $f_{\lambda}$ has the same inf and the same set of
minimizers as $f$ does, and that $f_{\lambda}$ has the same
symmetry properties as $f$ (that is, if $T:M\to M$ is an isometry
and $f(Tz)=f(z)$ for all $z$, then
$f_{\lambda}(Tz)=f_{\lambda}(z)$ for all $z$).

In Section 3 we assume that $f:M\to\mathbb{R}\cup\{+\infty\}$ is a
convex function and we study under what conditions on $M$ the
functions $f_{\lambda}$ are convex and $C^1$ smooth. It turns out
that some assumptions on the geometry of $M$ are necessary in
order that the $f_{\lambda}$ be convex and $C^1$ smooth (see
Example \ref{example sphere} below); in particular we must require
that the distance function $d:M\times M\to\mathbb{R}$,
$(x,y)\mapsto d(x,y)$, be {\em uniformly locally convex on bounded
sets near the diagonal} (see Definition \ref{definition of
uniformly locally convex near the diagonal} below). Under this
assumption we prove that the functions $f_{\lambda}$ are convex
and $C^1$ smooth on any given bounded subset $B$ of $M$, for all
$\lambda$ small enough. Moreover, if the distance $d$ is convex on
all of $M\times M$ then the $f_{\lambda}$ are convex and $C^1$
smooth on all of $M$ for all $\lambda>0$.

In Section 4 we study the question as to which manifolds satisfy
the above technical assumption that $d$ is uniformly locally
convex on bounded sets near the diagonal (resp. convex on $M\times
M$). First, we show that for every Riemannian manifold with
nonpositive sectional curvature, and with the property that the
convexity radius function of $M$ is strictly positive on bounded
sets (such is the case, for instance, of all complete
finite-dimensional Riemannian manifolds), the distance $d:M\times
M\to\mathbb{R}$ is uniformly locally convex on bounded sets near
the diagonal. Secondly, we note that, for every Cartan-Hadamard
manifold (that is, every simply connected complete Riemannian
manifold of nonpositive sectional curvature), the distance $d$ is
convex on all of $M\times M$. By combining these facts with the
results of Sections 2 and 3 we obtain, in the finite-dimensional
case (see Corollary \ref{main corollary} below):

{\em If $M$ is a complete finite-dimensional Riemannian manifold
with sectional curvature $K\leq 0$ and $f:M\to\mathbb{R}$ is a
convex function, then, for every bounded open convex set $U$ with
compact closure $\overline{U}$, there exists $\lambda_{0}>0$ such
that the functions $f_{\lambda}:M\to\mathbb{R}$ are convex and
$C^{1}$ smooth on $U$ for all $\lambda\in (0, \lambda_{0})$.
Moreover,
\begin{enumerate}
\item $f_{\lambda}$ converges to $f$ uniformly on $\overline{U}$.
\item $f_{\lambda}\leq f$ for all $\lambda>0$.
\item $f_{\lambda}$ has the same inf and the same set of minimizers
as $f$.
\item $f_{\lambda}$ has the same symmetries as $f$ (that is, if
$f$ is invariant with respect to an isometry $T:M\to M$, then so
is $f_{\lambda}$).
\end{enumerate}
}

And, in the case of a Cartan-Hadamard manifold (either
finite-dimensional  or infinite-dimensional, see Corollary
\ref{second main corollary} below):

{\em If $M$ is a Cartan-Hadamard manifold and
$f:M\to\mathbb{R}\cup\{+\infty\}$ is a lower-semicontinuous convex
function, then the functions $f_{\lambda}:M\to\mathbb{R}$ are
convex and $C^{1}$ smooth on all of $M$ for all $\lambda>0$.
Moreover,
\begin{enumerate}
\item If $f$ is uniformly continuous on bounded sets then $f_{\lambda}$ converges
to $f$ uniformly on bounded sets.
\item $f_{\lambda}\leq f$ for all $\lambda>0$.
\item $f_{\lambda}$ has the same inf and the same set of minimizers
as $f$.
\item $f_{\lambda}$ has the same symmetries as $f$.
\end{enumerate}
}

Finally, in Section 5 we consider some corollaries and
applications of the above results. We show that if $C$ is a closed
convex subset of a Cartan-Hadamard manifold then the distance
function to $C$, $x\mapsto d(x,C)=\inf\{d(x,y): y\in C\}$, is
$C^{1}$ smooth on $M\setminus C$ and the function $x\mapsto
d(x,C)^{2}$ is $C^1$ smooth and convex on all of $M$. We also note
that this result is not true for Riemannian manifolds of positive
curvature such as the $2$-sphere, and therefore the results of
Section 3 cannot be extended to manifolds of positive curvature.
Another consequence is that every closed convex subset of a
Cartan-Hadamard manifold can be approximated by $C^1$ smooth
convex bodies of $M$. Lastly, we note that if $M$ is a
Cartan-Hadamard manifold and $f:M\to\mathbb{R}\cup\{+\infty\}$ is
convex and lower-semicontinuous, then the function
$$
    u(t,x):=\inf_{y\in M}\{f(y)+\frac{1}{2t}d(x,y)^2\} \textrm{
    for } t>0, \,\, u(0,x)=f(x)
    $$
is the unique viscosity solution to the following Hamilton-Jacobi
partial differential equation:
    $$
    \begin{cases}
    \frac{\partial{u(t,x)}}{\partial t}+
    \frac{1}{2}\|\frac{\partial{u(t,x)}}{\partial x}\|^{2}_{x}=0 &  \\
    u(0,x)=f(x), &
    \end{cases}
    $$
where $u:[0,\infty)\times M\to\mathbb{R}$.

\medskip

\section{General properties}

\noindent Throughout the paper, for a function
$f:M\to\mathbb{R}\cup\{+\infty\}$, we
 define
    $$
    f_{\lambda}(x)=\inf_{y\in M}\{f(y)+\frac{1}{2\lambda}d(x,y)^{2}\}.
    $$
The following Proposition shows how, under certain conditions, the
inf defining $f_{\lambda}$ can be localized on a neighborhood of
the point $x$.
\begin{prop}\label{localization}[Localization]
Let $M$ be a Riemannian manifold,
$f:M\to\mathbb{R}\cup\{+\infty\}$ a function satisfying that
$f(x)\geq -\frac{c}{2}(1+d(x, x_{0})^{2})$ for some $c>0$,
$x_{0}\in M$. Let $x\in M$ be such that $f(x)<+\infty$. Then, for
all $\lambda\in (0,\frac{1}{2c})$ and for all $\rho>\bar{\rho}$,
where $$\bar{\rho}=\bar{\rho}(x, \lambda,
c):=\left(\frac{2f(x)+c(2 d(x,x_{0})^{2}+1)}{1-2\lambda
c}\right)^{1/2},$$ we have that
$$
f_{\lambda}(x)=\inf_{y\in
B(x,\rho)}\{f(y)+\frac{1}{2\lambda}d(x,y)^{2}\}.
$$
\end{prop}
\begin{proof}
Since $$d(y,x_{0})^{2}\leq\left( d(y,x)+ d(x,x_{0})\right)^{2}\leq
2\left( d(x,y)^{2}+d(x, x_{0})^{2}\right),$$ we have that
$$c\left( d(x,y)^{2}+d(x,
x_{0})^{2}\right)\geq\frac{c}{2}d(y,x_{0})^{2},$$ hence
$$
-\frac{c}{2}-\frac{c}{2}d(y,
x_{0})^{2}+\frac{1}{2\lambda}d(x,y)^{2}\geq\frac{1}{2\lambda}d(x,y)^{2}
-cd(x,y)^{2}-cd(x,x_{0})^{2}-\frac{c}{2},
$$
that is
$$-\frac{c}{2}\left(1+d(y,x_{0})^{2}\right)+\frac{1}{2\lambda}d(x,y)^{2}\geq
\left(\frac{1}{2\lambda}-c\right)d(x,y)^{2}-cd(x,x_{0})^{2}-\frac{c}{2}.
\eqno(1)
$$
Now, for any given $\eta>0$, if we set $$r=r(x,\lambda,c,\eta):=
\left(\lambda\frac{2f(x)+2\eta+c(2d(x,x_{0})^{2}+1)}{1-2\lambda
c}\right)^{1/2},$$ by using $(1)$ we obtain that, for every $y\in
M$ with $d(y,x)>r$,
\begin{eqnarray*}
& & f(y)+\frac{1}{2\lambda}d(x,y)^{2}\geq
-\frac{c}{2}\left(1+d(y,x_{0})^{2}\right)+\frac{1}{2\lambda}d(x,y)^{2}\geq\\
& &
\left(\frac{1}{2\lambda}-c\right)d(x,y)^{2}-cd(x,x_{0})^{2}-\frac{c}{2}\geq\\
& &\left(\frac{1}{2\lambda}-c\right)\left(\lambda\frac{2f(x)+
2\eta+c(2d(x,x_{0})^{2}+1)}{1-2\lambda
c}\right)-cd(x,x_{0})^{2}-\frac{c}{2}=\\
& & f(x)+\eta>f_{\lambda}(x),
\end{eqnarray*}
which implies that
$$
\inf_{d(y,x)>r}\{
f(y)+\frac{1}{2\lambda}d(x,y)^{2}\}>f_{\lambda}(x),
$$
hence
$$
f_{\lambda}(x)=\inf_{d(y,x)\leq r}\{
f(y)+\frac{1}{2\lambda}d(x,y)^{2}\}.
$$
Finally, since $$\lim_{\eta\to
0}r(x,\lambda,c,\eta)=\bar{\rho}(x,\lambda,c),$$ it is clear that
for every $\rho>\bar{\rho}$ we can find $\eta>0$ small enough so
that
$$\rho=\rho(x,\lambda,c)>r(x,\lambda,c,\eta)>\bar{\rho}(x,\lambda,c),$$
and therefore, from the above argument we deduce that
$$
f_{\lambda}(x)=\inf_{d(y,x)\leq \rho}\{
f(y)+\frac{1}{2\lambda}d(x,y)^{2}\}.
$$
\end{proof}

Next we state several interesting properties of this
inf-convolution operation, such as preservation of order and
symmetry properties of the original function. We put off studying
the conditions under which convexity is preserved until the next
section.
\begin{prop}
Let $M$ be a Riemannian manifold,
$f:M\to\mathbb{R}\cup\{+\infty\}$ a function. We have that:
\begin{enumerate}
\item $f_{\lambda}\leq f$ for all $\lambda>0$.

\item If $0<\lambda_{1}<\lambda_{2}$ then $f_{\lambda_{2}}\leq
f_{\lambda_{1}}$.

\item $\inf f_{\lambda}= \inf f$ and, moreover, if $f$ is lower
semicontinuous then every minimizer of $f_{\lambda}$ is a
minimizer of $f$, and conversely.

\item If $T$ is an isometry of $M$ onto $M$, and $f$ is invariant
under $T$ (that is, $f(Tz)=f(z)$ for all $z\in M$), then
$f_{\lambda}$ is also invariant under $T$, for all $\lambda>0$.
\end{enumerate}
\end{prop}
\begin{proof}
\noindent $(1)$ and $(2)$ are obvious.

\medskip

\noindent $(3)$ Note that
\begin{eqnarray*}
& &\inf_{x\in M}f_{\lambda}(x)=\inf_{x\in M}\inf_{y\in
M}\{f(y)+\frac{1}{2\lambda}d(x,y)^{2}\}=\\
& &\inf_{y\in M}\inf_{x\in
M}\{f(y)+\frac{1}{2\lambda}d(x,y)^{2}\}=\inf_{y\in M}f(y).
\end{eqnarray*}

\noindent $(4)$ We have that
\begin{eqnarray*}
& & f_{\lambda}(Tx)=\inf_{y\in
M}\{f(y)+\frac{1}{2\lambda}d(Tx,y)^{2}\}=\inf_{y\in
M}\{f(Ty)+\frac{1}{2\lambda}d(Tx,Ty)^{2}\}=\\
& &\inf_{y\in M}\{f(Ty)+\frac{1}{2\lambda}d(x,y)^{2}\}= \inf_{y\in
M}\{f(y)+\frac{1}{2\lambda}d(x,y)^{2}\}=f_{\lambda}(x).
\end{eqnarray*}
\end{proof}

\medskip

Now we apply Proposition \ref{localization} to show that, under
natural continuity assumptions on $f$, the regularizations
$f_{\lambda}$ converge to the original function $f$ as $\lambda$
goes to $0$.
\begin{prop}\label{convergence}[Convergence]
Let $M$ be a Riemannian manifold,
$f:M\to\mathbb{R}\cup\{+\infty\}$ a function satisfying that
$f(x)\geq -\frac{c}{2}(1+d(x, x_{0})^{2})$ for some $c>0$,
$x_{0}\in M$, and consider
$$f_{\lambda}(x)=\inf_{y\in
M}\{f(y)+\frac{1}{2\lambda}d(x,y)^{2}\}$$ for $0<\lambda<1/2c$.
\begin{enumerate}
\item Assume that $f$ is uniformly continuous on bounded subsets
of $M$. Then $\lim_{\lambda\to 0}f_{\lambda}=f$ uniformly on each
bounded subset of $M$.

\item Assume that $f$ is continuous on $M$. Then $\lim_{\lambda\to
0}f_{\lambda}=f$ uniformly on compact subsets of $M$.

\item Assume that $f$ is uniformly continuous and bounded on all
of $M$. Then $\lim_{\lambda\to 0}f_{\lambda}=f$ uniformly on $M$.

\item In general (that is, with no continuity assumptions on $f$)
we have that $\lim_{\lambda\to 0}f_{\lambda}(x)=f(x)$ for every
$x\in M$ with $f(x)<+\infty$.
\end{enumerate}
\end{prop}
\begin{proof}
\noindent (1) According to Proposition \ref{localization}, for
every $x\in M$, $\lambda\in (0, 1/2c)$,
$\rho>\bar{\rho}(x,\lambda,c)$, we have that $$
f_{\lambda}(x)=\inf_{d(x,y)\leq\rho}\{f(y)+\frac{1}{2\lambda}d(x,y)^{2}\}.$$
Fix $R>0$. As is easily shown, a uniformly continuous function on
a Riemannian manifold is bounded on bounded sets, hence we can
find $k>0$ so that $|f(x)|\leq k$ for all $x\in B(x_{0},2R)$.

For any given $\varepsilon>0$, by uniform continuity of $f$, there
exists $\delta>0$ such that if $y,x\in B(x_{0}, 2R)$ and
$d(x,y)\leq\delta$ then $|f(x)-f(y)|\leq\varepsilon/3$. We can
assume $\delta<R$. Now, since
$$
\lim_{\lambda\to
0^{+}}\left(\lambda\frac{2k+c(2R^{2}+1)}{1-2c\lambda}\right)^{1/2}=0,
$$
there is $\lambda_{\varepsilon}>0$ such that if
$0<\lambda<\lambda_{\varepsilon}$ then
$$
0<\rho(x,\lambda,c)\leq\left(\lambda\frac{2k+c(2R^{2}+1)}{1-2c\lambda}\right)^{1/2}<\delta
$$
for all $x\in B(x_{0}, R)$, and therefore
$$
f_{\lambda}(x)=\inf_{d(y,x)\leq\delta}\{f(y)+\frac{1}{2\lambda}d(y,x)^{2}\}
$$
for all $x\in B(x_{0}, R)$, $\lambda\in
(0,\lambda_{\varepsilon})$. But, since $f_{\lambda}\leq f$ for all
$\lambda$, this really means that
$$
f_{\lambda}(x)=\inf_{y\in A_{x}}\{
f(y)+\frac{1}{2\lambda}d(y,x)^{2}\},
$$
where $$A_{x}:=\{y\in B(x,\delta) :
f(y)+\frac{1}{2\lambda}d(y,x)^{2}\leq f(x)\}.$$ By the definition
of inf, we can take $y_{x}\in A_{x}$ such that
$$
f_{\lambda}(x)+\frac{\varepsilon}{3}\geq
f(y_{x})+\frac{1}{2\lambda}d(y_{x},x)^{2}.
$$
Then, bearing in mind that $y_{x}\in A_{x}\subseteq
B(x,\delta)\subseteq B(x_{0}, 2R)$ when $x\in B(x_{0},R)$, we get
\begin{eqnarray*}
& & |f(x)-f_{\lambda}(x)|=f(x)-f_{\lambda}(x)\leq f(x)-f(y_{x})-
\frac{1}{2\lambda}d(y_{x},x)^{2}+\frac{\varepsilon}{3}\leq\\
& &|f(x)-f(y_{x})|+|f(x)-f(y_{x})|+\frac{\varepsilon}{3}\leq
\frac{\varepsilon}{3}+\frac{\varepsilon}{3}+\frac{\varepsilon}{3}=\varepsilon
\end{eqnarray*}
for all $x\in B(x_{0}, R)$, $\lambda\in (0,
\lambda_{\varepsilon})$. This shows that $\lim_{\lambda\to
0^{+}}f_{\lambda}=f$ uniformly on $B(x_{0}, R)$.

\medskip

\noindent (2) Let $K$ be a compact subset of $M$. By compactness,
it is easily seen that for every $\varepsilon>0$ there exists
$\delta>0$ such that if $x\in K$, $y\in M$, and $d(x,y)\leq\delta$
then $|f(x)-f(y)|\leq\varepsilon/3$. One can now repeat the above
argument with the precaution of always taking $x\in K\subset
B(x_{0}, R)$, $y\in M$, $d(x,y)\leq\delta$.

\medskip

\noindent (3) Choose $k>0$ such that $|f(x)|\leq k$ for all $x\in
M$. Let us first observe that the inf defining $f_{\lambda}$ can
be restricted to the set $\{y\in M: d(y,x)\leq
2\sqrt{k\lambda}\}=B(x, 2\sqrt{k\lambda})$. Indeed, if
$d(y,x)>2\sqrt{k\lambda}$ then
$$f(y)+\frac{1}{2\lambda}d(x,y)^{2}>-k+2k=k\geq f(x)\geq
f_{\lambda}(x).$$

Next, for any given $\varepsilon>0$, the uniform continuity of $f$
provides us with $\delta>0$ so that $|f(y)-f(x)|\leq\varepsilon/3$
whenever $d(x,y)\leq\delta$. On the other hand, since
$\lim_{\lambda\to 0^{+}}2\sqrt{k\lambda}=0$, there exists
$\lambda_{\varepsilon}>0$ such that $0<2\sqrt{k\lambda}<\delta$
for $0<\lambda<\lambda_{\varepsilon}$. Then, for any $x\in M$,
$\lambda\in (0,\lambda_{\varepsilon})$ the inf defining
$f_{\lambda}(x)$ can be restricted to the set $A_{x}:=\{y\in
B(x,\delta): f(y)+\frac{1}{2\lambda}d(x,y)^{2}\leq f(x)\}$. Now,
the same estimations as in $(1)$ above show that
$|f_{\lambda}(x)-f(x)|\leq\varepsilon$.

\medskip

\noindent (4) is very easy.
\end{proof}

\medskip

\section{Regularization of convex functions}

In order to see that the operation $f\to f_{\lambda}$ preserves
convexity we will need to use the following Lemmas.

\begin{lem}\label{the partial inf of a jointly convex is convex}
Let $M$ be a Riemannian manifold, and $F:M\times
M\to\mathbb{R}\cup\{+\infty\}$ a convex function (where $M\times
M$ is endowed with its natural product Riemannian metric). Assume
either that $M$ has the property that every two points can be
connected by a geodesic in $M$, or else that $F$ is continuous and
$M$ is complete. Then, the function $\psi:M\to\mathbb{R}$ defined
by
$$\psi(x)=\inf_{y\in M}F(x,y)$$ is also convex.
\end{lem}
\begin{proof}
Let $\gamma:I\to M$ be a geodesic. We have to see that the
function $t\mapsto \psi(\gamma(t))$ is convex, that is, for any
$t_{0}, t_{1}\in I$, and for any $s\in [0,1]$,
$\psi(\gamma(st_{1}+(1-s)t_{0}))\leq
s\psi(\gamma(t_{1}))+(1-s)\psi(\gamma(t_{0}))$. We may assume (up
to an affine change of parameters) that $t_{0}=0$ and $t_{1}=1$,
so we have to show that $\psi(\gamma(t))\leq
t\psi(x_{1})+(1-t)\psi(x_{0})$, where $x_{0}:=\gamma(0)$ and
$x_{1}:=\gamma(1)$.

Fix an arbitrary $t\in [0,1]$. For any $\varepsilon>0$, by the
definition of $\psi$, we can pick $y_{0}, y_{1}\in M$ so that
$$F(x_{1}, y_{1})< \psi(x_{1})+\varepsilon, \,\,\, \textrm{ and }\,\, F(x_{0},
y_{0})< \psi(x_{0})+\varepsilon. \eqno(*)$$ Let $\sigma:J\to M$ be
a geodesic connecting $y_{0}$ and $y_{1}$ (if such $\sigma$ does
not exist then we can use continuity of $F$ and Ekeland's
approximate Hopf-Rinow type theorem to get points $\bar{y_{0}}$
and $\bar{y_{1}}$ close enough to $y_{0}$ and $y_{1}$ so that
$(*)$ remains true when replacing $y_{j}$ with $\bar{y_{j}}$, and
a geodesic $\bar{\sigma}$ connecting $\bar{y_{0}}$ to
$\bar{y_{1}}$; the rest of the argument applies without changes).
We can also assume that $J=[0,1]$, $y_{0}=\sigma(0)$,
$y_{1}=\sigma(1)$.

It is clear that, because $\gamma$ and $\sigma$ are geodesics in
$M$, the path $t\mapsto (\gamma(t),\sigma(t))$ is a geodesic
joining the points $(x_{0}, y_{0})$ and $(x_{1}, y_{1})$ in the
product manifold $M\times M$. Now, since $t\mapsto
F(\gamma(t),\sigma(t))$ is convex, we have that
\begin{eqnarray*}
& &\psi(\gamma(t))=\inf_{y\in M}F(\gamma(t), y)\leq
F(\gamma(t),\sigma(t))\leq t F(x_{1},y_{1})+(1-t)F(x_{0},
y_{0})\leq\\ & & t(\psi(x_{1})+\varepsilon)+
(1-t)(\psi(x_{0})+\varepsilon)=
t\psi(x_{1})+(1-t)\psi(x_{0})+\varepsilon
\end{eqnarray*}
and this holds for every $\varepsilon>0$, hence we can conclude
that $\psi(\gamma(t))\leq t\psi(x_{1})+(1-s)\psi(x_{0})$.
\end{proof}

\begin{lem}\label{convex functions are quadratically minored}
Let $M$ be a Riemannian manifold with the property that any two
points of $M$ can be joined by a minimizing geodesic, and let
$f:M\to\mathbb{R}\cup\{+\infty\}$ be a convex function. Then, for
every $x_{0}\in M$ there exists a number $c\geq 0$ such that
$f(x)\geq -\frac{c}{2}\left(1+d(x,x_{0})^{2}\right)$ for all $x\in
M$.
\end{lem}
\begin{proof}
Choose $\zeta\in D^{-}f(x_{0})$. For a given $x\in M$, let
$\gamma$ be a minimizing geodesic connecting $x_{0}$ to $x$, say
$\gamma(t)=\exp_{x_{0}}(tv)$, $t\in [0, d(x,x_{0})]$, for some
$v\in TM_{x_{0}}$ with $\|v\|_{x_{0}}=1$. Since $f$ is convex we
have $f(\exp_{x_{0}}(tv))-f(x_{0})\geq\langle\zeta,
tv\rangle_{x_{0}}$ for every $t\in [0, d(x_{0},x)]$, and in
particular
\begin{eqnarray*}
& & f(x)-f(x_{0})\geq \langle\zeta, d(x_{0},x)v\rangle_{x_{0}}\geq
-\|\zeta\|_{x_{0}}d(x_{0},x)\|v\|_{x_{0}}=\\
&
&-\|\zeta\|_{x_{0}}d(x_{0},x)\geq-\|\zeta\|_{x_{0}}\left(1+d(x_{0},x)^{2}\right),
\end{eqnarray*}
and therefore
$$
f(x)\geq f(x_{0})-\|\zeta\|_{x_{0}}\left(1+d(x_{0},x)^{2}\right)\geq
-\frac{c}{2}\left( 1+d(x_{0},x)^{2}\right)
$$
for all $x\in M$ if we put
$c=2\left(\|\zeta\|_{x_{0}}+|f(x_{0})|\right)$.
\end{proof}

It will be also useful to recall that every convex function
$f:M\to\mathbb{R}$ which is locally bounded is continuous (in fact
locally Lipschitz); a proof of this statement can be found in
\cite[Proposition 5.2]{AFL2}.

\medskip

In order that $f_{\lambda}$ is convex whenever $f$ is, we will
have to require that the distance function $d:M\times
M\to\mathbb{R}$ is convex on a band around the diagonal of
$M\times M$. More precisely, we will use the following.

\begin{defn}\label{definition of uniformly locally convex near the diagonal}
{\em Let $M$ be a Riemannian manifold. We say that the distance
function $d:M\times M\to\mathbb{R}$ is uniformly locally convex on
bounded sets near the diagonal if, for every bounded subset $B$ of
$M$ there exists $r>0$ such that $d$ is convex on $B(x,r)\times
B(x,r)$, and the set $B(x,r)$ is convex in $M$, for all $x\in
B)$.}
\end{defn}
Examples of manifolds satisfying this definition are the cylinder
$x^{2}+y^{2}=1$ in $\mathbb{R}^{3}$, the Poincar\'{e} half-plane, or
the subsets of $\mathbb{R}^{3}$ defined by $z=1/(x^{2}+y^{2})$ or
$z=xy$. In general, every complete finite-dimensional Riemannian
manifold of nonpositive sectional curvature meets this condition,
as we will show in the next section.

\begin{prop}\label{convexity preservation}
Let $M$ be a Riemannian manifold with the property that any two
points of $M$ can be joined by a geodesic, and let
$f:M\to\mathbb{R}\cup\{+\infty\}$ be a lower-semicontinuous convex
function.
\begin{enumerate}
\item  Assume that $f$ is bounded on
bounded sets and that the distance function $d:M\times
M\to\mathbb{R}$ is uniformly locally convex on bounded sets near
the diagonal. Then, for every bounded subset $B$ of $M$ there
exists $\lambda_{0}>0$ such that $f_{\lambda}$ is convex on $B$
for all $\lambda\in (0,\lambda_{0})$.

\item Assume that the distance function $d:M\times
M\to\mathbb{R}$ is convex on all of $M\times M$. Then $f_{\lambda}$
is convex on $M$ for every
$\lambda>0$.
\end{enumerate}
\end{prop}
\begin{proof}
\noindent $(1)$ We may well assume $B=B(x_{0}, R)$ for some
$x_{0}\in M$, $R>0$. Let $r>0$ be small enough so that the
function $(x,y)\mapsto d(y,x)$ is convex on $B(x, 2r)\times
B(x,2r)$ and $B(x,2r)$ is convex for every $x\in B(x_{0},R)$. Let
$k$ be a bound for $f$ on $B(x_{0}, 2R)$. We may assume $2r<R$. We
have that, for every $z\in B(x_{0},2R)$,
$$\bar{\rho}(z,\lambda,k):=\left(\lambda\frac{2f(z)+k(1+2R^{2})}{1-2\lambda
k}\right)^{1/2}\leq\left(\lambda\frac{2k+k(1+2R^{2})}{1-2\lambda
k}\right)^{1/2}\to 0 \,\, \text{ as } \, \lambda\to 0^{+}$$ hence we
can choose $\lambda_{0}>0$ small enough so that
$\bar{\rho}(z,\lambda_{0},k)<r$ for all $z\in B(x_{0}, 2R)$ and
therefore, according to Proposition \ref{localization}, we have
that, for every $\lambda\in
(0,\lambda_{0})$,$$f_{\lambda}(z)=\inf_{y\in
B(z,r)}\{f(y)+\frac{1}{2\lambda}d(z,y)^{2}\} = \inf_{y\in
B(x,2r)}\{f(y)+\frac{1}{2\lambda}d(z,y)^{2}\}$$whenever $z\in
B(x,r)$, $x\in B(x_{0},R)$. Note that $B(x,2r)$, as a convex subset
of $M$, still has the property that any two of its points can be
joined by a geodesic inside $B(x, 2r)$. Now, assuming
$0<\lambda<\lambda_{0}$, and fixing $x\in B(x_{0},R)$, because the
function $F(z,y):=f(y)+\frac{1}{2\lambda}d(z,y)^{2}$ is jointly
convex on $B(x,2r)\times B(x, 2r)$, we deduce from Lemma \ref{the
partial inf of a jointly convex is convex} that $z\mapsto
f_{\lambda}(z)=\inf_{y\in B(x,2r)}F(z,y)$ is convex on $B(x,r)$, for
all $0<\lambda<\lambda_{0}$. Since $x\in B(x_{0},R)$ is arbitrary
this implies that $f_{\lambda}$ is convex on $B(x_{0},R)$, for every
$\lambda\in (0, \lambda_{0})$.

\medskip

\noindent $(2)$ Here we can use Lemma \ref{the partial inf of a
jointly convex is convex} on all of $M\times M$ with no need to
localize the inf defining $f_{\lambda}(x)$, so it follows that
$f_{\lambda}$ is convex for all $\lambda>0$.
\end{proof}

\begin{rem}
{\em Note that in Case (2) of the above proposition we do not
require continuity of the function $f$, so $f$ is permitted to
take the value $+\infty$ at some points; in particular we are
allowed to take $f$ to be the indicator function of a closed
convex subset $C$ of $M$, that is
    $$
    \delta_{C}(x)=
  \begin{cases}
    0 & \text{ if } x\in C, \\
    +\infty & \text{ otherwise }.
  \end{cases}
    $$
}
\end{rem}

\begin{rem}
{\em An examination of the above proof and the statement of Lemma
\ref{the partial inf of a jointly convex is convex} reveals that,
if one assumes that $f$ is continuous and $M$ is complete, it is
not necessary to require that every two points of $M$ can be
connected by a minimizing geodesic in $M$.}
\end{rem}

The following Proposition shows that the functions $f_{\lambda}$
are superdifferentiable at a point $x$ if $d$ is differentiable on
a suitable ball around $x$. We refer the reader to \cite{AFL2} for
the properties of viscosity subdifferentials on Riemannian
manifolds; here we will only make use of the very definition of
the subdifferential and the superdifferential sets of
$f:M\to\mathbb[-\infty, \infty]$, namely,
    $$
    D^{-}f(x)=\{d\varphi(x): \varphi\in C^{1}(M,\mathbb{R}), \,\,
    f-\varphi \textrm{ attains a local minimum at } x\},
    $$
and $$
    D^{+}f(x)=\{d\psi(x): \psi\in C^{1}(M,\mathbb{R}), \,\,
    f-\psi \textrm{ attains a local maximum at } x\},
    $$
the fact that $f$ is differentiable at $x$ if and only if
$D^{-}f(x)\neq\emptyset\neq D^{+}f(x)$ (in which case $\{df(x)\}=
D^{-}f(x)= D^{+}f(x)$), and that a convex function
$f:M\to\mathbb{R}$ is everywhere subdifferentiable \cite[Theorem
5.3]{AFL2}
\begin{prop}\label{the regularizations are always superdifferentiable}
Suppose that the inf defining $f_{\lambda}(x)$
can be restricted to a ball $B_{x}=B(x,r_{x})$ of radius $r_{x}$
small enough so that the function $y\mapsto d(y,x)^{2}$ is
differentiable on $B_{x}$, and that this inf is attained at a point
$y_{x}\in B_{x}$. Then $f_{\lambda}$ is superdifferentiable at $x$,
and $$\frac{1}{\lambda}d(x, y_{x})\frac{\partial}{\partial x} d(x,
y_{x})\in D^{+}f_{\lambda}(x).$$
\end{prop}
\begin{proof}
We have that
\begin{eqnarray*}
& &f_{\lambda}(z)-f_{\lambda}(x)\leq
f(y_{x})+\frac{1}{2\lambda}d(z,y_{x})^{2}-f(y_{x})-\frac{1}{2\lambda}d(x,y_{x})^{2}=\\
&
&\frac{1}{2\lambda}d(z,y_{x})^{2}-\frac{1}{2\lambda}d(x,y_{x})^{2},
\end{eqnarray*}
so $$f_{\lambda}(z)-\frac{1}{2\lambda}d(z,y_{x})^{2}\leq
f_{\lambda}(x)-\frac{1}{2\lambda}d(x,y_{x})^{2}$$ for every $z\in
B_{x}$, that is, $f_{\lambda}-\frac{1}{2\lambda}d(\cdot, y_{x})^{2}$
attains a local maximum at $x$, hence $f_{\lambda}$ is
superdifferentiable at $x$, with $d\left( \frac{1}{2\lambda}d(\cdot,
y_{x})^{2}\right)(x)\in D^{+}f_{\lambda}(x)$.
\end{proof}
Next we show that convex differentiable functions are
automatically of class $C^1$. In this proof we will make use of
the parallel transport of vectors along geodesics. Recall that,
for a given curve $\gamma: I\to M$, numbers  $t_{0}, t_{1}\in I$,
and a vector $V_{0}\in TM_{\gamma(t_{0})}$, there exists a unique
parallel vector field $V(t)$ along $\gamma(t)$ such that
$V(t_{0})=V_{0}$. Moreover, the mapping defined by $V_{0}\mapsto
V(t_{1})$ is a linear isometry between the tangent spaces
$TM_{\gamma(t_{0})}$ and $TM_{\gamma(t_{1})}$, for each $t_{1}\in
I$. In the case when $\gamma$ is a minimizing geodesic and
$\gamma(t_{0})=x$, $\gamma(t_{1})=y$, we will denote this mapping
by $L_{xy}$, and we call it the parallel transport from $TM_{x}$
to $TM_{y}$ along the geodesic $\gamma$. See \cite{Klingenberg}
for general reference on these topics. The parallel transport
allows us to measure the length of the ``difference" between
vectors (or forms) which are in different tangent spaces (or in
duals of tangent spaces, that is, fibers of the cotangent bundle),
and do so in a natural way. Indeed, let $\gamma$ be a minimizing
geodesic connecting two points $x, y\in M$, say $\gamma(t_{0})=x,
\gamma(t_{1})=y$. Take vectors $v\in TM_{x}$, $w\in TM_{y}$. Then
we can define the distance between $v$ and $w$ as the number
    $$
    \|v-L_{yx}(w)\|_{x}=
    \|w-L_{xy}(v)\|_{y}
    $$
(this equality holds because $L_{xy}$ is a linear isometry between
the two tangent spaces, with inverse $L_{yx}$). Since the spaces
$T^{*}M_{x}$ and $TM_{x}$ are isometrically identified by the
formula $v=\langle v, \cdot\rangle$, we can obviously use the same
method to measure distances between forms $\zeta\in T^{*}M_{x}$
and $\eta\in T^{*}M_{y}$ lying on different fibers of the
cotangent bundle.

It is also well known that the mapping $y\mapsto L_{xy}$ is well
defined and continuous on a neighborhood of each $x\in M$, in the
following sense: if $(x_{n})$ converges to $x$ in $M$ then
$\exp_{x_{n}}(L_{xx_{n}}(v))$ converges to $\exp_{x}(v)$ uniformly
on the set $\{v\in TM_{x}: \|v\|_{x}\leq \delta\}$ for some
$\delta>0$ (a fact which we use at the end of the proof of the
following lemma).

\begin{lem}\label{differentiable convex functions are C1}
Let $M$ be a Riemannian manifold, and let $f:M\to\mathbb{R}$ be a
differentiable {\em convex} function. Then $f$ is of class $C^1$ on
$M$.
\end{lem}
\begin{proof}
Assume that $f$ is not $C^1$, then there are $\varepsilon>0$, a
point $x\in M$ and a sequence $(x_{n})\subset M$ converging to $x$
such that
$$\|L_{x_{n}x}[df(x_{n})]-df(x)\|_{x}>2\varepsilon$$ for all
$n\in\mathbb{N}$. Therefore, for every $n\in\mathbb{N}$ we can pick
$h_{n}\in TM_{x}$ with $\|h_{n}\|_{x}=1$ such that
$$\langle
L_{x_{n}x}[df(x_{n})]-df(x),h_{n}\rangle_{x}>2\varepsilon\,\,
\textrm{ for all }\,n\in\mathbb{N}.$$
Since $f$ is differentiable at
$x$, there exists $\delta>0$ so that
$$f(\exp_{x}(tv))-f(x)-\langle
df(x), tv\rangle_{x}\leq \varepsilon t$$ for all $v\in TM_{x}$ with
$\|v_{x}\|=1$ and $|t|\leq\delta$. On the other hand, by convexity
of $f$, we have $$\langle df(x_{n}), tw\rangle_{x_{n}}\leq
f(\exp_{x_{n}}(tw))-f(x_{n})$$ for all $w\in TM_{x_{n}}$ with
$\|w\|_{x_{n}}=1$ and $|t|\leq\delta$. By combining these
inequalities we get
\begin{eqnarray*}
& & 2\varepsilon\delta\leq \langle
L_{x_{n}x}[df(x_{n})]-df(x),h_{n}\rangle_{x}\delta=\langle
L_{x_{n}x}[df(x_{n})], \delta h_{n}\rangle_{x}-\langle df(x),
\delta h_{n}\rangle_{x}=\\& &\langle  df(x_{n}), \delta
L_{xx_{n}}h_{n}\rangle_{x_{n}}-\langle df(x), \delta
h_{n}\rangle_{x}\leq \\ & & f(\exp_{x_{n}}(\delta
L_{xx_{n}}h_{n}))-f(x_{n})+f(x)-f(\exp_{x}(\delta
h_{n}))+\varepsilon\delta=\\& &f(\exp_{x_{n}}(\delta
L_{xx_{n}}h_{n}))-f(\exp_{x}(\delta
h_{n}))+f(x)-f(x_{n})+\varepsilon\delta \to 0+0+\varepsilon\delta
\end{eqnarray*}
(by continuity of $f$, $\exp$ and the parallel translation $L$), so
we get that $2\varepsilon\delta\leq\varepsilon\delta$, which is not
possible.
\end{proof}
Now we can prove that, under the same assumptions as in Proposition
\ref{convexity preservation}, if $f$ is convex then $f_{\lambda}$ is
of class $C^1$ for $\lambda>0$ small enough.
\begin{thm}\label{regularization of convex functions}
Let $M$ be a Riemannian manifold and let
$f:M\to\mathbb{R}\cup\{+\infty\}$ be a lower semicontinuous and
convex function. Assume that every two points of $M$ can be
connected by a minimizing geodesic in $M$.
\begin{enumerate}
\item Suppose that $f$ is bounded on
bounded sets and that the distance function $d:M\times
M\to\mathbb{R}$ is uniformly locally convex on bounded sets near
the diagonal. Then, for every bounded open convex subset $B$ of
$M$ there exists $\lambda_{0}>0$ such that $f_{\lambda}$ is a
$C^1$ smooth convex function on $B$, for all $\lambda\in
(0,\lambda_{0})$.

\item Suppose that the distance function $d:M\times
M\to\mathbb{R}$ is convex on all of $M\times M$. Then
$f_{\lambda}$ is a $C^1$ smooth convex function on $M$, for every
$\lambda>0$.
\end{enumerate}
\end{thm}
\begin{proof}
\noindent $(1)$ We can give an almost self-contained proof of this
in the finite-dimensional case, so let us first assume that
$\textrm{dim} M<+\infty$. We may also assume $B=B(x_{0},R)$ for
some $x_{0}\in M$, $R>0$. Since the index of injectivity $x\mapsto
i(x)$ is a continuous positive function, it is bounded below by a
positive number on the compact subset $\overline{B}(x_{0},R)$ of
$M$. This implies that there exists $r_{0}>0$ such that the
function $y\mapsto d(x,y)^{2}$ is $C^1$ smooth on $B(x,r_{0})$ for
every $x\in B(x_{0},R)$. We can obviously assume that $r<r_{0}$
and repeat the argument of the proof of Proposition \ref{convexity
preservation} to get a $\lambda_{0}>0$ such that $f_{\lambda}$ is
convex on $B(x_{0},R)$ for all $\lambda\in (0,\lambda_{0})$ and,
moreover, that, for every $x\in B(x_{0}, R)$, the inf defining
$f_{\lambda}(x)$ can be restricted to the ball
$\overline{B}(x,r)$, which is contained in $B(x, r_{0})$ (so that,
in particular, $y\mapsto d(y,x)^{2}$ is $C^1$ smooth on $B(x,r)$).
Besides, this inf is attained, because the involved functions are
continuous and the ball $\overline{B}(x,r)$ is compact. According
to Proposition \ref{the regularizations are always
superdifferentiable}, we then get that $f_{\lambda}$ is
superdifferentiable at $x$, for every $x\in B(x_{0},R)$,
$\lambda\in (0,\lambda_{0})$.

On the other hand, since $f_{\lambda}$ is convex on $B(x_{0},R)$, we
know that $f_{\lambda}$ is subdifferentiable on $B(x_{0},R)$ (see
\cite[Theorem 5.3]{AFL2}). That is, $f_{\lambda}$ is both
subdifferentiable and superdifferentiable at each point of
$B(x_{0},R)$, hence $f_{\lambda}$ is differentiable on $B(x_{0},R)$
(see \cite[Proposition 4.6]{AFL2}). Since $f_{\lambda}$ is convex,
Lemma \ref{differentiable convex functions are C1} allows to
conclude that $f_{\lambda}$ is $C^1$ smooth on $B(x_{0},R)$ for each
$\lambda\in (0,\lambda_{0})$.

\medskip

Let us now consider the case when $\textrm{dim} M=+\infty$. Since
$M$ has the property that any two of its points can be connected
by a minimizing geodesic in $M$, and $f_{\lambda}$, being convex
on $B$ for all $\lambda\in (0,\lambda_{0})$, satisfies
$D^{-}f_{\lambda}(x)\neq\emptyset$ for all $x\in B$, we can apply
Theorem 11 of \cite{AF1} to get that $f_{\lambda}$ is
differentiable at every point $x\in B$, hence (by Lemma
\ref{differentiable convex functions are C1}) of class $C^1$ on
$B$.

\medskip

\noindent $(2)$ As in case $(1)$, let us first give a
self-contained proof for the finite-dimensional case. If the
distance function $d:M\times M\to\mathbb{R}$ is convex on all of
$M\times M$ then $y\mapsto d(y,x)$ is convex on $M$ for all $x\in
M$, which implies that the cut locus of $x$ is empty for every
$x\in M$ and that the function $y\mapsto d(y,x)^{2}$ is
differentiable on all of $M$ for every $x\in M$. On the other
hand, we claim that the inf defining $f_{\lambda}(x)$ is attained
for every $x\in M$. Indeed, fix $x\in M$ with $f(x)<+\infty$. From
the proof of Lemma \ref{convex functions are quadratically
minored} we know that there exists $c=c_{x}\geq 0$ such that
$f(y)\geq-c d(x,y)$ for every $y\in M$. Then we have
$$
f(y)+\frac{1}{2\lambda}d(x,y)^{2}\geq-c d(x,y)+
\frac{1}{2\lambda}d(x,y)^{2}\to+\infty
$$
if $d(x,y)\to+\infty$, so there exists $R>0$ large enough so that
if $d(x,y)\geq R$ then
    $$
    f(y)+\frac{1}{2\lambda}d(x,y)^{2}\geq f(x)\geq f_{\lambda}(x),
    $$
hence $$f_{\lambda}(x)=\inf_{y\in
\overline{B}(x,R)}\{f(y)+\frac{1}{2\lambda}d(x,y)^{2}\},
$$
and now it is clear that this inf is attained because
$f+\frac{1}{2\lambda}d(\cdot, x)^{2}$ is lower semicontinuous and
$\overline{B}(x,R)$ is compact.

Therefore, according to Proposition \ref{the regularizations are
always superdifferentiable}, $f_{\lambda}$ is superdifferentiable.
Because $f_{\lambda}$ is convex, this means that $f_{\lambda}$ is
diferentiable, hence $C^1$ smooth on $M$, for all $\lambda>0$.

\medskip

In the infinite-dimensional case we only have to bear in mind
that, according to Proposition \ref{convexity preservation}(2),
$f_{\lambda}$ is convex on all of $M$, so the same proof as in
case $(1)$ applies.
\end{proof}
\begin{rem}
{\em If one assumes that $f$ is continuous and $M$ is complete, it
is not necessary to require that every two points of $M$ can be
connected by a minimizing geodesic in $M$.}
\end{rem}

\medskip

\section{Which manifolds do the above results apply to?}

Let us say a few words about the manifolds satisfying the
assumptions of Theorem \ref{regularization of convex functions}.
The following Theorem is a restatement of \cite[Theorem IX. 4.3,
p. 257]{Lang}

\begin{thm}\label{Theorem from Lang's book}
Let $M$ be a Riemannian manifold with seminegative sectional
curvature $K\leq 0$, and $U$ a convex open set. Let $\beta_{1},
\beta_{2}$ be disjoint geodesics in $U$, defined on the same
interval. Let $\alpha_{t}:[a,b]\to U$ be the unique geodesic
joining $\beta_{1}(t)$ with $\beta_{2}(t)$, and let $\ell
(t)=\textrm{length}(\alpha_{t})$, that is,
$\ell(t)=d(\beta_{1}(t), \beta_{2}(t))$. Then $\ell''(t)\geq 0$
for all $t$, and in particular $\ell(t)$ is a convex function.
\end{thm}

From this Theorem it is immediate to deduce that the above results
on regularization of convex functions apply to manifolds of
seminegative sectional curvature, as we next see.

\begin{cor}
Let $M$ be a Riemannian manifold with sectional curvature $K\leq
0$.
\begin{enumerate}
\item  Suppose that $M$ has a convexity radius function which is
strictly positive on bounded subsets of $M$ (such is the case, for
instance, of a complete finite-dimensional Riemannian manifold
$M$). Then the distance function $d$ is uniformly locally convex
on bounded sets near the diagonal of $M\times M$.
\item Suppose that $M$ is simply connected (which, together with
the curvature assumption, amounts to saying that $M$ is a
Cartan-Hadamard manifold). Then the distance function $d$ is
convex on all of $M$.
\end{enumerate}
\end{cor}
\begin{proof}
$(1)$ Let $B$ be a compact subset of $M$. Since the convexity
radius function $x\mapsto c(x)$ is bounded below on $B$ by a
number $r>0$, we have that the open ball $B(x,r)$ is convex for
every $x\in B$. Therefore, for every $x\in B$ and for every pair
of disjoint geodesic segments $\beta_{1}, \beta_{2}:I:=(a,b)\to M$
contained in $B(x,r)$, Theorem \ref{Theorem from Lang's book}
tells us that the function $t\mapsto \ell(t):=d(\beta_{1}(t),
\beta_{2}(t))$ is convex. If $\beta_{1}, \beta_{2}$ are not
disjoint and neither of them is constant (in which case the result
would be trivial) then we can only have the equality
$d(\beta_{1}(t),\beta_{2}(t))=0$ for a unique $t=t_{0}$, at which
point the function $\ell(t)$ attains an absolute minimum, and
Theorem \ref{Theorem from Lang's book} shows that $\ell(t)$ is
convex on $(a,t_{0})$ and on $(t_{0},b)$. But a real function
which is convex on $(a,t_{0})$ and on $(t_{0},b)$, and which
attains its minimum at $t_{0}$, must in fact be convex on all of
$I=(a,b)$.

This proves that the distance function $d$ is convex on
$B(x,r)\times B(x,r)$, for every $x\in B$, which in turn means
that $d$ is uniformly locally convex near the diagonal.

\medskip

$(2)$ In a Cartan-Hadamard manifold $M$ every ball is convex, and
two distinct geodesics in $M$ can intersect in only one point (see
\cite[p. 259-261]{Lang}), so the above argument applies globally.
\end{proof}
\begin{rem}
{\em The assumption on curvature is necessary in order that $d$ be
uniformly locally convex near the diagonal: it is easy to see
that, for many disjoint nonconstant geodesic segments $\beta_{1}$
and $\beta_{2}$ in the sphere $S^{2}$ (take for instance two
parallel meridians near the equator), the function $t\to
d(\beta_{1}(t), \beta_{2}(t))$ is {\em not} convex. Furthermore,
as we will see in the next section, an important consequence of
Theorem \ref{regularization of convex functions} fails in the
sphere $S^{2}$, so the assumption on the jointly convexity of the
distance function $d:M\times M\to\mathbb{R}$ near the diagonal
seems to be much more than a mere technical requirement and is
probably a necessary condition for the functions $f_{\lambda}$ to
be convex whenever $f$ is.}
\end{rem}

\medskip

We conclude with a Corollary that sums up what the results we have
shown tell us in the case of a Riemannian manifold of nonpositive
curvature.

When $M$ is a complete finite-dimensional Riemannian manifold of
nonpositive curvature we have the following result. Recall that a
convex function $f$ on a finite-dimensional Riemannian manifold
$M$ that only takes finite values is automatically continuous (see
\cite{Bangert}).

\begin{cor}\label{main corollary}
Let $M$ be a complete finite-dimensional Riemannian manifold with
sectional curvature $K\leq 0$. Let $f:M\to\mathbb{R}$ be a convex
function. Then, for every bounded open convex set $U$ with compact
closure $\overline{U}$, there exists $\lambda_{0}>0$ such that the
functions $f_{\lambda}:M\to\mathbb{R}$, defined by
$$
f_{\lambda}(x)=\inf_{y\in M}\{f(y)+\frac{1}{2\lambda}d(x,y)^{2}\},
$$
are convex and $C^{1}$ smooth on $U$ for all $\lambda\in (0,
\lambda_{0})$. Moreover,
\begin{enumerate}
\item $f_{\lambda}$ converges to $f$ uniformly on $\overline{U}$.
\item $f_{\lambda}\leq f$ for all $\lambda>0$.
\item $f_{\lambda}$ has the same inf and the same set of minimizers
as $f$.
\item $f_{\lambda}$ has the same symmetries as $f$ (that is, if
$f$ is invariant with respect to an isometry $T:M\to M$, then so
is $f_{\lambda}$).
\end{enumerate}
\end{cor}

Note that in this result we do not allow $f$ to take infinite
values. We are able to deal with functions
$f:M\to\mathbb{R}\cup\{+\infty\}$ when we furthermore assume that
$M$ is a Cartan-Hadamard manifold (that is, a simply connected
Riemannian manifold of nonpositive curvature), either finite or
infinite-dimensional. Also recall that Cartan-Hadamard manifolds
enjoy the property that every two points can be connected by a
minimizing geodesic (see \cite{Lang}).
\begin{cor}\label{second main corollary}
Let $M$ be a Cartan-Hadamard manifold (either finite-dimensional
or infinite-dimensional). Let $f:M\to\mathbb{R}\cup\{+\infty\}$ be
a lower-semicontinuous convex function. Then the functions
$f_{\lambda}:M\to\mathbb{R}$, defined by
$$
f_{\lambda}(x)=\inf_{y\in M}\{f(y)+\frac{1}{2\lambda}d(x,y)^{2}\},
$$
are convex and $C^{1}$ smooth on all of $M$ for all $\lambda>0$.
Moreover,
\begin{enumerate}
\item $f_{\lambda}\leq f$ for all $\lambda>0$.
\item $f_{\lambda}$ has the same inf and the same set of minimizers
as $f$.
\item $f_{\lambda}$ has the same symmetries as $f$ does.
\item If $f$ is uniformly continuous on bounded sets $f_{\lambda}$ converges to $f$
uniformly on bounded sets.
\end{enumerate}
\end{cor}

\medskip

\section{Some applications}

If $X$ is a Hilbert space (or more generally a reflexive Banach
space), it is well known that for every closed convex subset $C$
of $X$ the distance function to $C$, that is, $x\mapsto d(x, C)$
is convex and $C^1$ smooth away from $C$ (even though $C$ might
not have a smooth boundary), and, as a consequence, every such $C$
can be approximated by $C^1$ smooth convex bodies. We next show
how the results proved above allow us to extend these two theorems
to the class of Cartan-Hadamard manifolds (either finite or
infinite-dimensional), and we also note that this result
completely fails in the sphere $S^2$: there are closed convex sets
$C$ of arbitrarily small diameter in $S^2$ such that $x\mapsto
d(x,C)$ is not convex on any neighborhood of $C$.

\begin{cor}\label{diferentiability of the distance to a convex compact}
Let $C$ be a closed convex subset of a Cartan-Hadamard manifold.
Then the distance function to $C$, $x\mapsto d(x,C)=\inf\{d(x,y):
y\in C\}$ is $C^{1}$ smooth on $M\setminus C$ and, moreover, the
function $x\mapsto d(x,C)^{2}$ is $C^1$ smooth and convex on all
of $M$.
\end{cor}
\begin{proof}
Define $f:M\to\mathbb{R}\cup\{+\infty\}$ by
    $$
    f(x)=
  \begin{cases}
    0 & \text{ if } x\in C, \\
    +\infty & \text{ otherwise }.
  \end{cases}
    $$
The function $f$ is lower semicontinuous and convex on $M$.
According to Theorem \ref{regularization of convex functions}, the
function $f_{\lambda}:M\to\mathbb{R}$,
    $$
    f_{\lambda}(x)=\inf\{f(y)+\frac{1}{2\lambda}d(x,y)^{2}\}=
    \inf\{\frac{1}{2\lambda}d(x,y)^{2}: y\in C\}=\frac{1}{2\lambda}d(x,C)^{2},
    $$
is $C^1$ smooth and convex on $M$ for all $\lambda>0$. By taking
$\lambda=1/2$ we get that the squared distance function to $C$ is
$C^1$ smooth and convex on $M$.
\end{proof}

\begin{defn}
{\em We say that a subset $C$ of a Riemannian manifold $M$ is a
$C^1$ smooth {\em convex body} of $M$ provided $C$ is closed,
convex, has nonempty interior, and $\partial C$ is a
one-codimensional $C^1$ smooth submanifold of $M$.}
\end{defn}

\begin{cor}\label{Smooth approximation of convex bodies}
Let $C$ be a closed convex subset of a Cartan-Hadamard manifold,
and let $U$ be an open subset of $M$ with $d(C, M\setminus U)>0$.
Then there exists a $C^1$ smooth convex body $D$ of $M$ such that
$C\subset D\subset U$.
\end{cor}
\begin{proof}
Since $d(C, \partial U)>0$ we can take $r=\frac{1}{2}d(C,\partial
U)$ and define $D=\{x\in M: d(x, C)\leq r\}$. It is clear that
$C\subset D\subset U$, and $D$ happens to be a $C^1$ smooth convex
body because $x\mapsto d(x,C)^{2}$ is $C^1$ smooth and convex, and
the derivative of the function $d(\cdot, C)^{2}$ is nonzero at
every point $x\in M\setminus C$ (a convex function has a null
derivative only at the points, if any, where it attains its
minimum).
\end{proof}

\medskip

The following example shows that Theorem \ref{regularization of
convex functions} and the above Corollaries are false in general
if we do not require that the manifold $M$ has nonpositive
sectional curvature.

\begin{ex}\label{example sphere}
{\em Let $M$ be the sphere $x^{2}+y^{2}+z^{2}=1$ in
$\mathbb{R}^{2}$ endowed with its usual Riemannian metric. Let $C$
be a closed geodesic segment of diameter less than a number
$\varepsilon$ with $0<\varepsilon<1$. It is easy to see that the
function $d(\cdot, C)^{2}$ defined on $M$ by
    $$
    d(x,C)^{2}=\inf_{y\in C} d(x,y)^{2},
    $$
is not convex on any open neighborhood of $C$, that is, Corollary
\ref{diferentiability of the distance to a convex compact} fails
in $M$ (hence so does Theorem \ref{regularization of convex
functions}).}
\end{ex}

\medskip

Finally, it should be noted that there is a strong link between
the regularization method we have just presented and the following
Hamilton-Jacobi partial differential equation:
    $$
  (*)\begin{cases}
    \frac{\partial{u(t,x)}}{\partial t}+
    \frac{1}{2}\|\frac{\partial{u(t,x)}}{\partial x}\|^{2}_{x}=0 &  \\
    u(0,x)=f(x), &
  \end{cases}
    $$
where $u:[0,\infty)\times M\to\mathbb{R}$,
$f:M\to\mathbb{R}\cup\{+\infty\}$. If we assume that $M$ is a
finite-dimensional Cartan-Hadamard manifold and $f$ is convex and
lower-semicontinuous, then the function
    $$
    u(t,x):=\inf_{y\in M}\{f(y)+\frac{1}{2t}d(x,y)^2\} \textrm{
    for } t>0, \,\, u(0,x)=f(x)
    $$
is the unique viscosity solution of $(*)$ (see \cite{AFL2} for the
definition of viscosity solution to Hamilton-Jacobi equations on
Riemannian manifolds). This is not very difficult to show
directly. Alternatively, one can prove that Theorem 3.6 and
Section 7.2 of \cite{Barles} remain true when $\mathbb{R}^{n}$ is
replaced with a finite-dimensional Cartan-Hadamard manifold.



\begin{thebibliography}{}

\bibitem{AF1} D. Azagra and J. Ferrera, {\em Proximal calculus on
Riemannian manifolds}, preprint, 2004.

\bibitem{AFL2} D. Azagra, J. Ferrera, F. L\'{o}pez-Mesas, {\em
Nonsmooth analysis and Hamilton-Jacobi equations on Riemannian
manifolds}, J. Funct. Anal. 220 (2005) no. 2, 304-361.

\bibitem{Bangert}
V. Bangert, {\em Analytische Eigenschaften konvexer Funktionen auf
Riemannschen Mannigfaltigkeiten}, J. Reine Angew. Math. 307/308
(1979), 309--324.

\bibitem{Bangert2}
V. Bangert, {\em \"{U}ber die Approximation von lokal konvexen
Mengen}, Manuscripta Math. 25 (1978), no. 4, 397--420.

\bibitem{Bangert3}
V. Bangert, {\em Totally convex sets in complete Riemannian
manifolds}, J. Differential Geom. 16 (1981), no. 2, 333--345.

\bibitem{Barles}
G. Barles, {\em Solutions de viscosit\'{e} des \'{e}quations de
Hamilton-Jacobi}, Math\'{e}matiques \& Applications 17,
Springer-Verlag, Berlin, 1994.

\bibitem{C1}
M. Cepedello Boiso, {\em Approximation of Lipschitz functions by
$\Delta$-convex functions in Banach spaces}, Israel J. Math. 106
(1998), 269--284.

\bibitem{C2}
M. Cepedello Boiso, {\em On regularization in superreflexive Banach
spaces by infimal convolution formulas}, Studia Math. 129 (1998),
no. 3, 265--284.

\bibitem{CheegerGromoll}
J. Cheeger, and D. Gromoll, {\em On the structure of complete
manifolds of nonnegative curvature}, Ann. of Math. 96 (1972),
413--443.

\bibitem{Greene1}
R. E. Greene, and K. Shiohama, {\em Convex functions on complete
noncompact manifolds: topological structure}, Invent. Math. 63
(1981), no. 1, 129--157.

\bibitem{Greene2}
R. E. Greene, and K. Shiohama, {\em Convex functions on complete
noncompact manifolds: differentiable structure}, Ann. Sci. \'{E}cole
Norm. Sup. (4) 14 (1981), no. 4, 357--367 (1982).

\bibitem{Greene5}
R. E. Greene, and H. Wu, {\em On the subharmonicity and
plurisubharmonicity of geodesically convex functions}, Indiana
Univ. Math. J. 22 (1972/73), 641--653.

\bibitem{Greene3}
R. E. Greene, and H. Wu, {\em $C\sp{\infty }$ convex functions and
manifolds of positive curvature}, Acta Math. 137 (1976), no. 3-4,
209--245.

\bibitem{Greene4}
R. E. Greene, and H. Wu, {\em $C\sp{\infty }$ approximations of
convex, subharmonic, and plurisubharmonic functions}, Ann. Sci.
\'{E}cole Norm. Sup. (4) 12 (1979), no. 1, 47--84.

\bibitem{GromollMeyer}
D. Gromoll, and W. Meyer, {\em On complete open manifolds of
positive curvature}, Ann. of Math. 90 (1969) 75--90.

\bibitem{Klingenberg}
W. Klingenberg, {\em Riemannian geometry}, de Gruyter Studies in
Mathematics, 1. Walter de Gruyter \& Co., Berlin-New York, 1982.

\bibitem{Lang}
S. Lang, {\em Fundamentals of Differential Geometry}, Graduate
Texts in Mathematics, 191. Springer-Verlag, New York, 1999.

\bibitem{Lasry-Lions}
J.-M. Lasry, and P.-L. Lions,  {\em A remark on regularization in
Hilbert spaces.} Israel J. Math. 55 (1986), no. 3, 257--266.

\bibitem{Udriste}
C. Udriste, {\em Convex functions and optimization methods on
Riemannian manifolds}, Mathematics and its Applications, 297. Kluwer
Academic Publishers, Dordrecht, 1994.

\end{thebibliography}
\end{document}